
\magnification=\magstep1

\font\BBF=msbm10
\def\Bbf#1{\hbox{\BBF #1}}
\font\SCR=cmsy10
\def\Scr#1{\hbox{\SCR #1}}

\def\P{\Bbf P}

\def\Z{\Bbf Z}
\def\A{\Bbf A}

\def\O{\Scr O}

\def\qed
{\hskip 10pt \hbox{
\vrule height 7.5pt depth -0.1pt \vrule height 7.53pt depth -7.2pt width 7.3pt
\hskip -7.5pt \vrule height 0.3pt depth 0pt width 7.6pt \vrule height 7.5pt
depth -0.1pt
}}

\def\wt{\widetilde}

\def\inv{\sp{-1}}

\def\dim{\hbox{\rm dim }}

\def\Sing{\hbox{\rm Sing\hskip 2pt}}

\def\aa{\alpha}

\def\ww{\omega}

\def\GG{\Gamma}

\def\smalln{\smallskip\noindent}
\def\medn{\medskip\noindent}
\def\bign{\bigskip\noindent}
\def\parn{\par\noindent}

\def\lra{\longrightarrow}

\def\surj{\to \hskip -7pt \to}

\def\hs{\hskip 5pt}

\def\varmaprightsp#1#2{\smash{\mathop{\hbox to #1 {\rightarrowfill}}
\limits\sp{#2}}}
\def\varmaprightsb#1#2{\smash{\mathop{\hbox to #1 {\rightarrowfill}}
\limits\sb{#2}}}

\def\hookdownarrow
{{}\sp\cap \hskip -3.703pt \lower 2pt \hbox{$\downarrow$}}
\def\hookuparrow
{\lower 2pt\hbox{${}\sb{\cup}$}\hskip -3.7pt \lower -1pt\hbox{$\uparrow$}}

\def\setbar{\hs ; \hs}

\def\set#1#2{\{\hs{#1}\setbar{#2}\}}

\def\sethd#1#2#3
{\Bigl\{ \hs{#1}\setbar
{\matrix{\hbox{#2} \hfill\cr\hbox{#3} \hfill }} \Bigr\} }
\def\setht#1#2#3#4
{\biggl\{ \hs{#1}\setbar
{\matrix{\hbox{#2} \hfill\cr\hbox{#3} \hfill\cr\hbox{#4} \hfill}}
\biggr\} }

\def\locus#1{\{ #1 \}}

\def\varvarmatrix#1#2#3#4
{\def\normalbaselines{\baselineskip#1\lineskip3pt\lineskiplimit3pt}
\vbox{\vskip #2 \hbox{\hfill$\matrix{#3}$\hfill}\vskip #3}}

\def\varmatrix#1#2
{\def\normalbaselines{\baselineskip#1\lineskip3pt\lineskiplimit3pt}
\matrix{#2}}

\def\bigcases#1
{\biggl\{\, \vcenter{\normalbaselines{\mathsurround=0pt}
\ialign{$##\hfil$&\quad##\hfil\crcr#1\crcr}}\biggr.}

\def\st{\subset}
\def\sm{\setminus}

\def\pione{\pi\sb1}

\def\vol#1{{\bf #1 }}

\def\and{\hbox{and}}
\def\where{\hbox{where}}

\def\afcite{\ }

\def\Pt{\P\sp 2}
\def\hc{X, Y, Z}
\def\F{\Scr F}
\def\L{\Scr{L}}

\def\W{\Scr{W}}

\def\pqi{\phi\sb q \inv}

\def\ArtalBartolo{[1]}
\def\Libgober{[2]}
\def\Nemethi{[3]}
\def\OkaS{[4]}
\def\OkaPre{[5]}
\def\ShimadaF{[6]}
\def\ShimadaW{[7]}
\def\Tokunaga{[8]}
\def\ZariskiE{[9]}
\def\ZariskiP{[10]}
\def\ZariskiD{[11]}

\def\ShimadaFPart#1{[6; #1]}
\def\ShimadaWPart#1{[7; #1]}

\def\sbinfty{\sb{\infty}}

\def\spno{\sp{n-1}}

%
%
%

\centerline{\bf A note on Zariski pairs}
\bigskip
\centerline{Ichiro Shimada}
\smallskip
\centerline{Max-Planck-Institut f\"ur Mathematik, Bonn}

\bign
{\bf \S 0. Introduction}
\medskip
In \ArtalBartolo,
Artal Bartolo defined the notion of {\sl Zariski pairs} as follows:
\medskip
{\bf Definition.}
\hs
A couple of complex reduced projective
plane curves $C\sb 1$ and $C\sb 2$ of a same degree
is said to make a  Zariski pair,
if there exist tubular neighborhoods
$T(C\sb i) \st \Pt $ of $C\sb i$ for $i=1,2$
such that $(T(C\sb 1), C\sb 1)$ and $(T(C\sb 2 ), C\sb 2)$
are diffeomorphic,
while the pairs $(\Pt, C\sb 1)$ and $(\Pt, C\sb 2)$ are  not homeomorphic;
that is, the singularities of $C\sb 1 $ and $C\sb 2$
are topologically equivalent,
but the embeddings of $C\sb 1$ and $C\sb 2$
into $\Pt$ are not topologically equivalent.
\medskip
The first example of Zariski pair was discovered and studied
by Zariski in \ZariskiE\afcite and \ZariskiD.
He showed that
there exist projective plane curves $C\sb 1 $ and $C\sb 2$ of degree $6$
with $6$ cusps and no other singularities
such that $\pione (\Pt \sm C\sb 1)$ and $\pione (\Pt \sm C\sb 2)$
are not isomorphic.
Indeed, the placement of the $6$ cusps  on the sextic curve  has a
crucial effect on
the fundamental group of the complement.
Let $C\sb 1$ be a sextic curve defined by an equation $f\sp 2 + g \sp 3= 0$,
where $f$ and $g$ are general homogeneous polynomials of degree $3$ and $2$,
respectively.
Then $C\sb 1$ has $6$ cusps lying on a conic defined by $g=0$.
In \ZariskiE, it was shown
that  $\pione (\Pt\sm C\sb 1)$ is isomorphic to the free product
$\Z/(2) * \Z/(3)$
of cyclic groups of order $2$ and $3$.
On the other hand,
in \ZariskiD, it was proved that there exists  a sextic curve $C\sb 2$
with $6$ cusps which are not  lying on any conic, and
that the fundamental group $\pione (\Pt \sm C\sb 2)$ is cyclic of order $6$.
In \OkaS, Oka gave an explicit defining equation of $C\sb 2$.
In \ArtalBartolo,
Artal Bartolo presented  a simple way to construct
$(C\sb 1, C\sb 2)$
from a cubic curve $C$
by means of a Kummer covering of $\P\sp 2$
of exponent $2$
branched along three lines tangent to $C$ at its points of inflection.
\medskip
Except for this example,
very few Zariski pairs are known (\ArtalBartolo, \Tokunaga).
In \OkaPre, and independently in \ShimadaW,
infinite series of Zariski pairs
have been constructed from the above example of Zariski by means of
covering tricks of the plane.
\medskip
In this paper,
we present a method to construct Zariski pairs,
which yields
two infinite series of new examples of Zariski pairs as special cases.
\medskip
A germ of curve singularity is called of type $(p, q)$
if it is locally defined by $x\sp p + y\sp q=0$.
\medskip
{\it Series I.}
\hs
This series consists of pairs $(C\sb 1 (q), C\sb 2 (q))$
of curves of degree $3q$, where
$q$ runs through the set of integers $\ge 2$
prime to $3$.
Each of $C\sb 1 (q)$ and $C\sb 2 (q)$
has $3q$ singular points of type $(3,q)$
and no other singularities.
The fundamental group $\pione (\P\sp 2 \sm C\sb 1 (q))$
is non-abelian,
while
$\pione (\P\sp 2 \sm C\sb 2 (q))$
is abelian.
When $q=2$,
this example is nothing but  the classical one
of the sextic curves due to Zariski.
\medskip
{\it Series II.} \hs
This series consists of pairs $(D\sb 1 (q), D\sb 2 (q))$
of curves of degree $4q$, where
$q$ runs through the set of  odd integers $> 2$.
Each of $D\sb 1 (q)$ and $D\sb 2 (q)$
has $8q$ singular points of type $(2,q)$ -- that is,
rational double  points of type $A\sb{q-1}$ --
and no other singularities.
The fundamental group $\pione (\P\sp 2 \sm D\sb 1 (q))$
is non-abelian,
while
$\pione (\P\sp 2 \sm D\sb 2 (q))$
is abelian.
\medskip
Our method is a generalization of Artal Bartolo's method
for re-constructiong the classical example of Zariski
to higher dimensions and arbitrary exponents of the Kummer covering.
Indeed, when $q=2$ in Series I,
our construction coincides with his.
\medskip
Instead of  the computation
of the first Betti number
of the cyclic branched covering
of $\P\sp 2$,
which was employed in \ArtalBartolo,
we use the fundamental groups of the complements
in order to distinguish two embeddings of curves in $\Pt$.
For the calculation of the fundamental groups,
we use \ShimadaFPart{Theorem 1}
and a result of \Nemethi\afcite and \ShimadaW.
\medskip
{\sl Acknowledgment.}
\hs
Part of this work was done during the author's stay
at Institute of Mathematics in Hanoi
and Max-Planck-Institut f\"ur Mathematik  in Bonn.
The author thanks to people at these institutes
for their warm hospitality.
He also thanks
to Professors M.\ Oka  and H.\ Tokunaga
for stimulating discussions.
\bign
{\bf \S 1. A method of constructing Zariski pairs}
\medskip
{\bf  1.1. Non-abelian members.}
\hs
Let $p$ and $q$ be integers $\ge 2$
prime to each other.
We choose homogeneous polynomials
$f \in H\sp 0 (\Pt, \O(pk))$ and $g \in H\sp 0 (\Pt, \O (qk))$,
where $k$ is an integer $\ge 1$.
Suppose that $f$ and $g$ are generally chosen.
Consider the projective plane curve
$$
C\sb{p,q, k} \quad : \quad f\sp q + g\sp p = 0
$$
of degree $pqk$ (cf.\hs  \Libgober).
It is easy to see that the singular locus of this curve
consists of $pqk\sp 2 $ points of type $(p, q)$.
In \ShimadaWPart{Example (3) in \S 0},
the following is shown.
\medn
{\bf Proposition 1.}
\hs
{\sl
The fundamental group $\pione (\Pt\sm C\sb{p, q, k})$
is isomorphic to the group
$\langle \hs a, b, c \hs | \hs a\sp p = b\sp q = c,
\hs c\sp k = 1 \hs \rangle$.
In particular, it is non-abelian.
}
\medn
See also \Nemethi, in which the fundamental groups
of the complements of curves of this type are  calculated.
There the groups are presented in a different way.
\medskip
This curve $C\sb{p, q, k}$ will be a member $C\sb 1$ of a Zariski pair.
\medskip
{\bf 1.2. Abelian partners.}
\hs
We shall construct the other  member $C\sb 2 $
of the Zariski pair such that
$\pione (\Pt \sm C\sb 2)$ is abelian.
\medskip
Let $p$, $q$ and $k$ be integers as above.
We put $n = pk$.
Interchanging $p$ and $q$ if necessary,
we may assume that $n \ge 3$.
Let $S\sb 0 \st \P\sp{n-1}$ be a hypersurface of degree $n$
defined by
$F\sb 0 (X\sb 1, \dots, X\sb n) = 0$.
We consider a linear pencil of hypersurfaces
$$
S\sb t \quad : \quad F\sb 0 (X\sb 1, \dots, X\sb n) +
t \cdot X\sb 1 \cdots X\sb n =0,
$$
which is spanned by $S\sb 0$ and $S\sb{\infty} :=
\locus {X\sb 1 \cdots X\sb n =0}$.
We put $H\sb i = \locus{X\sb i = 0}$ $(i= 1, \dots, n)$.
We consider the morphism
$\phi\sb q : \P\sp {n-1} \to \P\sp{n-1}$
given by
$$
(Y\sb 1 : \dots : Y\sb n) \quad \longmapsto
\quad ( X\sb 1 : \dots : X\sb n ) \hs =\hs
( Y\sb 1 \sp q : \dots : Y\sb n \sp q),
$$
which is a covering of degree $q\spno$ branched along $S\sbinfty$.
\medn
{\bf Proposition 2.}
\hs
{\sl
Suppose that (1) every member $S\sb t$
is reduced,
and that (2) $S\sb 0$ contains none of the hyperplanes $H\sb i$.
Then
$\pione (\P\spno\sm \phi\sb q \inv (S\sb t))$ is abelian
for a general member $S\sb t$.
}
\medn
{\it Proof.}
\hs
Let $\P\sp 1$ be the $t$-line,
and we put $\A\sp 1 := \P\sp 1 \sm \locus {\infty}$.
Let $\W \st \P\spno \times \A\sp 1$
be the divisor defined by
$$
X\sb 1 \cdots X\sb n \cdot
(\hs  F\sb 0 (X\sb 1 , \dots , X\sb n) + t\cdot X\sb 1 \cdots X\sb n \hs ) = 0,
$$
which is the union of $S\sbinfty \times \A\sp 1 $
and the universal family of the
affine part $\set{ S\sb t }{t \in \A\sp 1}$
of the pencil.
For $t \in \A\sp 1$,
we denote by $W\sb t \st \P\spno$
the divisor obtained from the scheme theoretic intersection
$ (\locus{t} \times \P\spno ) \cap \W$,
which is equal with the divisor $S\sb t + S\sbinfty$.
\smallskip
First, we shall show that $\pione (\P\spno \sm W\sb t)$
is abelian for a general $t$.
Remark  that the assumption (2) implies that $S\sb t$ contains
none of $H\sb i$ unless $t=\infty$.
Combining this with the assumption (1),
we see that $W\sb t$ is reduced for all $t\in \A\sp 1$.
Hence, by \ShimadaFPart{Theorem 1},
the inclusion
$\P\spno \sm W\sb t  \hookrightarrow (\P\spno \times \A\sp 1 ) \sm \W$
induces an isomorphism on the fundamental groups
for a general $t$.
Therefore,
it is enough to show that
$\pione ((\P\spno \times \A\sp 1 )\sm \W)$
is abelian.
In order to prove this,
we consider the first projection
$$
p \quad :\quad
(\P\spno \times \A\sp 1 )\sm \W \quad \lra \quad \P\spno \sm S\sbinfty.
$$
Since $\set{ S\sb t }{t \in \P\sp 1}$
is a pencil whose base locus is contained in $S\sbinfty$,
there is a unique member $S\sb{t(P)}$ $(t(P)\not = \infty)$
containing $P$
for each point $P\in \P\spno \sm S\sbinfty$.
Therefore
$p\inv (P)$ is a punctured affine line $\A\sp 1 \sm\locus {t(P)}$
for every $P\in \P\spno \sm S\sbinfty$.
Consequently,  $p$ is a locally trivial fiber space.
Moreover,
$p$ has a section
$$
s \quad : \quad
\P\spno \sm S\sbinfty \quad\lra\quad (\P\spno \times \A\sp 1 )\sm \W,
$$
which is given by, for example,
$s(P) = (P, \hs t(P) +1)$.
Hence the homotopy exact sequence of $p$ splits.
Combining this with the fact that the image of the injection
$\pione (\A\sp 1 \sm \{ t(P) \}) \to \pione ((\P\spno \times \A\sp 1 )\sm\W)$
is contained in the center, we see that
$$
\pione ((\P\spno \times \A\sp 1 )\sm \W)
\quad \cong\quad
\pione (\P\spno \sm S\sbinfty )\hs
\times \hs \pione (\A\sp 1 \sm \locus {\hbox { a point }}).
$$
This shows that
$\pione ((\P\spno \times \A\sp 1 )\sm \W)$
is abelian.
\smallskip
Note that
$\phi\sb q : \P\spno \to \P\spno$
is \'etale over $\P\spno \sm W\sb t$ for every $t$.
Hence the natural homomorphism
$$
\phi\sb{q*} \quad :
\quad \pione (\P\spno \sm \phi\sb q \inv (W\sb t)) \quad \lra \quad
\pione (\P\spno \sm W\sb t)
$$
is injective.
This implies that
$\pione (\P\spno\sm \phi\sb q \inv (W\sb t))$
is abelian for a general $t$.
On the other hand,
since
$\P\spno \sm \phi\sb q \inv (W\sb t)$
is a Zariski open dense subset of $\P\spno \sm \phi\sb q \inv (S\sb t)$,
the inclusion
induces a surjective homomorphism
$$
\pione (\P\spno \sm \phi\sb q \inv (W\sb t)) \quad\surj\quad
\pione (\P\spno \sm \phi\sb q \inv (S\sb t)).
$$
Thus  $\pione (\P\spno \sm \phi\sb q \inv (S\sb t))$ is also
abelian for a general $t$.
\qed
\medn
{\bf Proposition 3.}
\hs
{\sl
Suppose the following;
(3) $S\sb 0 \cap H\sb i$
is a non-reduced divisor $p D\sb i$ of $H\sb i$ of multiplicity $p$,
where $D\sb i$ is a reduced divisor of $H\sb i$,
none of  whose irreducible components
is contained in $H\sb i \cap (\cup\sb{j\not = i} H\sb j)$,
and
(4) the singular locus of $S\sb t$ is of codimension $\ge 2$ in $S\sb t$
for a general $t$.
Then the general plane section
$\P\sp 2 \cap \phi \sb q \inv (S\sb t)$ of $\pqi (S\sb t)$
is a curve of degree $pqk$,
and its singular locus consists of $pqk\sp 2 $ points of type $(p, q)$.
}
\medn
{\it Proof.}
\hs
Note that the assumption (3)
implies that $S\sb t \cap H\sb i$
is also equal with $p D\sb i$ for $t\not = \infty$.
Let $P$ be a general point of any irreducible component of $D\sb i$,
and let $Q$ be a point such that $\phi\sb q (Q) = P$,
which lies on the hyperplane defined by $Y\sb i = 0$.
By the assumption (3),
$Q$ is not contained in any of the other hyperplanes
defined by $Y\sb j = 0$ $(j\not = i)$.
Hence there exist analytic local coordinate systems
$(w\sb 1 , \dots, w\sb{n-1})$
and $(z\sb 1 , \dots, z\sb{n-1})$
of $\P\spno$ with the origins $P$ and $Q$,
respectively,
such that
$H\sb i$ is given by $w\sb 1 =0$,
$\phi\sb q \inv (H\sb i)$
is given by $z\sb 1 =0$,
and $\phi\sb q$
is given by
$$
(z\sb 1 , \dots, z\sb{n-1} ) \quad \longmapsto
\quad (w\sb 1, \dots, w\sb{n-1}) \hs =
\hs  (z\sb 1 \sp q , z\sb 2 , \dots, z\sb{n-1}).
$$
Let $t \in \A\sp 1$ be general.
By the assumption (3),
the defining equation of $S\sb t$ at $P$ is of the form
$$
u(w)\cdot w\sb 1 + v(w\sb 2, \dots, w\sb{n-1}) \sp p = 0.
$$
By the assumption (4),
$S\sb t$ is non-singular at $P$,
because $P$ is a general point
of an irreducible component of $D\sb i$.
This implies that
$u(P) \not= 0$.
On the other hand,
the divisor $D\sb i$,
which is defined by $v(w\sb 2 , \dots, w\sb{n-1})=0$
on the hyperplane $H\sb i = \locus{w\sb 1 =0}$,
is non-singular at $P$,
because $D\sb i$ is reduced
by the assumption (3) and $P$ is general.
Hence we have
$$
{{\partial v}\over{\partial w\sb j}} (P) \not= 0
\quad\hbox{ at least for one $j \ge 2$.}
$$
The defining equation of $\phi\sb q\inv (S\sb t)$
is then of the form
$$
\wt u (z) \cdot z\sb 1 \sp q + v (z\sb 2 , \dots, z\sb{n-1} )\sp p = 0,
\quad\where\quad \wt u (Q) \not = 0.
$$
Then, it is easy to see that,
in terms of suitable analytic coordinates $(\wt z\sb 1 , \dots,
\wt z\sb{n-1})$
with the origin $Q$,
this equation can be written
as follows;
$$
\wt z\sb 1 \sp{\hskip 2pt q} \hs + \hs  \wt z \sb 2 \sp{\hskip 2pt p}
\hs = \hs 0.
$$
Thus, when we cut $\phi\sb q \inv (S\sb t)$
by a general $2$-dimensional plane passing through $Q$,
a germ of curve singularity of type $(p, q)$ appears at $Q$.
\smallskip
Since the degree of $D\sb i$ is $k = n/p$,
the inverse image $\phi\sb q \inv (D\sb i)$
is a reduced hypersurface of degree $qk$ in
the hyperplane defined by $Y\sb i=0$.
Moreover $\pqi (D\sb i)$ and $\pqi (D\sb j)$
have no common irreducible components
when $i\not = j$
because of the assumption (3).
Hence  the intersection  points of
$\phi\sb q \inv (\sum \sb{i=1} \sp n D\sb i)$
with a general plane $\Pt \st \P\sp n$
is $pqk\sp 2$ in number.
Moreover, $\Pt \cap \phi\sb q \inv (S\sb t)$
is non-singular outside of these intersection points,
because of the assumption (4).
\qed
\medskip
{\bf  1.3. Summary.}
\hs
Suppose that we have constructed a hypersurface $S\sb 0 \st \P\spno$
of degree $n\ge 3$ which satisfies the assumptions (1)-(4)
in Propositions 2  and 3.
Let $C\sb 2$ be a general plane section of $\pqi (S\sb t)$,
where $t$ is general.
Because of Zariski's hyperplane section theorem
\ZariskiP\afcite  and Propositions 1, 2 and 3,
we see that
the curve $C\sb 2$ has the same type of singularities as that of
$C\sb{ p, q, k}$,
but the fundamental group $\pione (\Pt \sm C\sb 2)$ is abelian.
Hence $(C\sb 1, C\sb 2)$ is a Zariski pair, with $C\sb 1 = C\sb{p, q, k}$.
\bign
{\bf \S 2.  Construction of Series I}
\medskip
We carry out the construction of the previous section with
$p=3$, $k=1$, $n=3$  and
$q$ an arbitrary integer $\ge 2$ prime to $3$.
\medskip
We fix a homogeneous coordinate system $(X:Y:Z)$
of $\Pt$, and put
$$
\varmatrix{20pt}{
&L\sb1  = \locus{X=0}, \quad L\sb2  = \locus{Y=0},
\quad L\sb 3  = \locus {Z=0}, \quad \and \cr
&R\sb 1 = (0:1:-1) \hs \in \hs L\sb 1,
\quad
R\sb 2 = (1:0:-1) \hs \in \hs L\sb 2.
}
$$
Let $\P\sb * (\GG (\Pt, \O(3)))$ be the space of all cubic curves on $\Pt$,
which is isomorphic to the projective space of dimension $9$,
and  let $\F \st \P\sb * (\GG (\Pt, \O(3)))$ be the family of
cubic curves $C$
which satisfy the following conditions;
\smalln
(a) $C$ intersects $L\sb 1$ at $R\sb 1$ with multiplicity $ \ge 3$,
\parn
(b) $C$ intersects $L\sb 2$ at $R\sb 2$ with multiplicity $\ge 3$, and
\parn
(c) $C$ intersects $L\sb 3$ at a point  with multiplicity $\ge 3$.
\smalln
(We consider that $C$ intersects a line $L\sb i$
with multiplicity $\infty$, if $L\sb i $ is contained in $C$.)
\medn
{\bf Proposition 4.}
\hs
{\sl
The family $\F$ consists of $3$ projective lines.
They meet at one point corresponding to
$C\sbinfty : = \locus {XYZ=0}$.
}
\medn
{\it Proof.}
\hs
Let $F(\hc)=0$ be the defining equation of a member $C$ of this family $\F$.
By the condition (a),
$F$ is of the form
$$
F(\hc) = A (Y+Z) \sp 3 + X \cdot G(\hc),
$$
where $A$ is a constant,
and $G(\hc)$ is a homogeneous polynomial of degree $2$.
By the condition (b), we have
$F(X, 0, Z) = A (Z+X)\sp 3$,
and hence we get
$$
G(\hc) = A (3Z\sp 2 + 3ZX + X\sp 2 ) + Y  \cdot H(\hc),
$$
where $H (\hc) $ is a homogeneous polynomial of degree $1$.
By the condition (c), we have
$F(X, Y, 0) = A (Y+ \aa X)\sp 3$
for some $\aa$. Then $\aa$ must be a cubic root of unity, and we get
$$
H(\hc) = 3A\aa \sp 2 X + 3A \aa Y + BZ,
$$
where $B$ is a constant.
Combining all of these, we get
$$
\eqalign{
&F(\hc) \cr
= &A(X\sp 3 + Y\sp 3 + Z\sp 3) + 3 A ( \aa\sp 2 X\sp 2 Y
+ \aa X Y \sp 2 +  Y\sp 2 Z +  Y Z \sp 2 +  Z\sp 2 X +  Z X \sp 2)
+ BXYZ\cr
=& A(X+Y+Z) \sp 3 +
3 A (\aa\sp 2 -1 ) X\sp 2 Y + 3 A (\aa -1) X Y\sp 2 + (B-6A) XYZ.
}
$$
This  curve  $C = \locus{F=0}$ intersects $L\sb 3$ at
$$
R\sb 3 \hs = \hs  R\sb 3 (\aa) \hs := \hs (1: -\aa: 0) \hs \in \hs L\sb 3
$$
with multiplicity $\ge 3$.
This means that the family $\F$
consists of three lines $\L (1)$, $\L(\ww)$ and $\L(\ww\sp 2)$
in the projective space  $\P\sb * (\GG (\Pt, \O(3)))$,
where $\ww= \exp(2\pi i / 3)$,
such that a general cubic
$C$ in $\L (\aa)$ intersects $L\sb 3$ at $R\sb 3 (\aa)$
with multiplicity $3$.
The ratio of the coefficients $t:= B/ A$
gives an affine coordinate on each line $\L (\aa)$.
The three lines $\L (1)$, $\L(\ww)$, $\L(\ww\sp 2)$
intersect at one point $t = \infty$
corresponding to the cubic $C\sb{\infty}=L\sb 1 + L\sb2 + L\sb 3$.
\qed
\medskip
Hence we get three pencils of cubic curves
$\set{ C(1)\sb t }{t \in \L (1)}$,
$\set{ C(\ww)\sb t }{t \in \L (\ww)}$, and
$\set{ C(\ww\sp 2 )\sb t }{t \in \L (\ww\sp 2)}$.
It is easy to check that these pencils
satisfy
the assumptions (2), (3) and (4)
in the previous section.
Note that  the pencil $\L (1)$
does not satisfy the assumption (1)
because $C(1)\sb 6 $ is a triple line.
However, the other two satisfy (1).
Indeed, if a cubic curve
$C$ in the family $\F$ is non-reduced,
then the conditions (a)-(c) imply
that it must be a triple line. Therefore the three points
$R\sb 1$, $ R\sb 2$ and $ R\sb 3 (\aa) $ are co-linear,
which is equivalent to $\aa=1$.
Consequently, $C$ must be a member of $\L (1)$.
\medskip
Now,
by using the pencil $\L (\ww)$
or $\L(\ww\sp 2)$, we complete  the construction of Series I.
\medskip
Note that,
if $C(1)\sb a$ is a non-singular member of $\L (1)$,
then
$\pione (\Pt \sm \phi\sb q \inv (C(1)\sb a))$
is isomorphic to the free product $\Z / (3) * \Z / (q)$.
Indeed, since
$C(1)\sb a$ is defined by
$$
(X+Y+Z)\sp 3 + (a-6)XYZ = 0,
$$
the pull-back
$\phi\sb q\inv (C(1)\sb a)$
is defined by
$$
(U\sp q + V\sp q + W\sp q)\sp 3 + (a-6) (UVW)\sp q =0,
$$
which is of the form $\wt f\sp{\hskip 2pt 3} + \wt g\sp{\hskip 2pt q} =0$.
The polynomials $\wt f$ and $\wt g$ are
not general by any means.
However,
since the type of singularities of $\phi\sb q \inv (C(1)\sb a )$
is the same as that of $C\sb{3,q,1}$,
we have an isomorphism
$\pione (\Pt\sm \phi\sb q \inv (C(1)\sb a)) \cong \pione
(\Pt \sm C\sb{3,q,1})$.
\bign
{\bf \S 3. Construction of Series II}
\medskip
It is enough to show the following:
\medn
{\bf Proposition 5.}
\hs
{\sl
The quartic surface
$$
S\sb 0 \quad:\quad
F\sb 0 (x\sb 1, x\sb 2, x\sb 3, x\sb 4) \hs :
= \hs (x\sb 1 \sp 2 + x\sb 2 \sp 2 )\sp 2 +
2 x\sb 3 x\sb 4 (x\sb 1 \sp 2 - x\sb 2 \sp 2 ) + x\sb 3 \sp 2 x\sb 4 \sp 2
\hs = \hs 0
$$
in $\P\sp 3 $ satisfies the assumptions (1)-(4) with $p=2$ and $k=2$.
}
\medn
{\it Proof.}
\hs
The assumptions (2) and (3) can be trivially checked.
To check the assumptions (1) and (4),
we put
$$
F\sb t  \hs  := \hs F\sb 0 + t\cdot x\sb 1 x\sb 2 x\sb 3 x\sb 4,
$$
and calculate the partial derivatives $\partial F\sb t / \partial x\sb i$
for $i=1, \dots, 4$.
Let $Q\sb t \st \P\sp 3$
be the quadric surface defined by
$$
 2 x\sb 1 \sp 2 - 2 x\sb 2 \sp 2 + 2 x\sb 3 x\sb 4 + t x\sb 1 x\sb 2 \hs
= \hs 0.
$$
It is easy to see that $Q\sb t$ is irreducible for all $t\not = \infty$.
It is also easy to see that $Q\sb t$ is the unique common
irreducible component
of the two cubic surfaces
$$
{{\partial F\sb t }\over{\partial x\sb 3 }} = 0,
 \qquad\and\qquad
{{\partial F\sb t }\over{\partial x\sb 4 }} = 0.
$$
Suppose that a surface $S\sb a = \locus{F\sb a =0}$
in this pencil contains a non-reduced
irreducible component $m T$ $(m \ge 2)$.
Then, both of $\partial F\sb a / \partial x\sb 3$ and
$\partial F\sb a / \partial x\sb 4$
must vanish on $T$.
Hence $T$ must coincide with $Q\sb a $, and we get $S\sb a = 2 Q\sb a$.
Comparing the defining equations of $S\sb a $ and $2 Q\sb a$,
we see that there are no such $a$.
Thus the assumption (1) is satisfied.
To check the assumption (4),
we remark that the condition $\dim\Sing S\sb t \le 0$
is an open condition for $t$.
Hence it is enough to prove, for example,
$\dim \Sing S\sb 2 =0$.
It is easy to show that $\Sing S\sb 2$
consists of four points
$(1:\pm\sqrt{-1}:0:0)$, $(0:0:0:1)$ and $(0:0:1:0)$.
\qed
\bigskip
\centerline{\bf References}
\medskip
\item{\ArtalBartolo}
 Artal Bartolo, E.:
{\sl Sur les couples de Zariski},
J.\ Alg.\ Geom. {\bf 3} (1994), 223 - 247
\item{\Libgober}
Libgober, A.:
{\sl Fundamental groups of the complements to plane singular curves},
Proc.\ Symp.\ in Pure Math.
{\bf 46} (1987), 29 - 45
\item{\Nemethi}
N\'emethi, A.:
{\sl On the fundamental group of the complement
of certain singular plane curves},
Math. Proc. Cambridge Philos. Soc.
{\bf 102} (1987), 453 - 457
\item{\OkaS}
Oka, M.:
{\sl Symmetric plane curves with nodes and cusps},
J.\ Math.\ Soc.\ Japan {\bf 44} (1992), 375 - 414
\item{\OkaPre}
Oka, M.: {\sl Two transformations of plane
curves and their fundamental groups},
preprint
\item{\ShimadaF}
Shimada, I.:
{\sl Fundamental groups of open algebraic varieties},
to appear in Topology
\item{\ShimadaW}
Shimada, I.:
{\sl A weighted version of Zariski's hyperplane section theorem
and fundamental groups of complements of plane curves},
preprint
\item{\Tokunaga}
Tokunaga, H.:
{\sl A remark on Bartolo's paper},
preprint
\item{\ZariskiE}
Zariski, O.:
{\sl On the problem of existence of
algebraic functions of two variables
possessing a given branch curve},
Amer. J. Math. {\bf 51} (1929), 305 - 328
\item{\ZariskiP}
Zariski, O.:
{\sl A theorem on the Poincar\'e group of an algebraic hypersurface},
Ann.\ Math. {\bf 38} (1937), 131 - 141
\item{\ZariskiD}
Zariski, O.:
{\sl The topological discriminant group of a Riemann surface
of genus $p$},
Amer.\ J.\ Math.\ \vol{59} (1937), 335 - 358
\bign
Max-Planck-Institut f\"ur Mathematik
\parn
Gottfried-Claren-Strasse 26
\parn
53225 Bonn, Germany
\parn
shimada@mpim-bonn.mpg.de

\end
----------
X-Sun-Data-Type: default
X-Sun-Data-Description: default
X-Sun-Data-Name: eprint
X-Sun-Charset: us-ascii
X-Sun-Content-Lines: 1053

\\
Title: A note on Zariski pairs
Author: Ichiro Shimada
TeX-Type: Plain TeX
Subj-Class: 14H30
Notes: 7 pages
\\
A couple of complex projective plane curves are said to make a Zariski pair
if they have the same degree and the same type of singularities, but their
embeddings in the projective plane are topologically different. In this paper,
we present a method to construct certain type of Zariski pairs, and make
infinite series of examples.
\\

\magnification=\magstep1

\font\BBF=msbm10
\def\Bbf#1{\hbox{\BBF #1}}
\font\SCR=cmsy10
\def\Scr#1{\hbox{\SCR #1}}

\def\P{\Bbf P}

\def\Z{\Bbf Z}
\def\A{\Bbf A}

\def\O{\Scr O}

\def\qed
{\hskip 10pt \hbox{
\vrule height 7.5pt depth -0.1pt \vrule height 7.53pt depth -7.2pt width 7.3pt
\hskip -7.5pt \vrule height 0.3pt depth 0pt width 7.6pt \vrule height 7.5pt
depth -0.1pt
}}

\def\wt{\widetilde}

\def\inv{\sp{-1}}

\def\dim{\hbox{\rm dim }}

\def\Sing{\hbox{\rm Sing\hskip 2pt}}

\def\aa{\alpha}

\def\ww{\omega}

\def\GG{\Gamma}

\def\smalln{\smallskip\noindent}
\def\medn{\medskip\noindent}
\def\bign{\bigskip\noindent}
\def\parn{\par\noindent}

\def\lra{\longrightarrow}

\def\surj{\to \hskip -7pt \to}

\def\hs{\hskip 5pt}

\def\varmaprightsp#1#2{\smash{\mathop{\hbox to #1 {\rightarrowfill}}
\limits\sp{#2}}}
\def\varmaprightsb#1#2{\smash{\mathop{\hbox to #1 {\rightarrowfill}}
\limits\sb{#2}}}

\def\hookdownarrow
{{}\sp\cap \hskip -3.703pt \lower 2pt \hbox{$\downarrow$}}
\def\hookuparrow
{\lower 2pt\hbox{${}\sb{\cup}$}\hskip -3.7pt \lower -1pt\hbox{$\uparrow$}}

\def\setbar{\hs ; \hs}

\def\set#1#2{\{\hs{#1}\setbar{#2}\}}

\def\sethd#1#2#3
{\Bigl\{ \hs{#1}\setbar
{\matrix{\hbox{#2} \hfill\cr\hbox{#3} \hfill }} \Bigr\} }
\def\setht#1#2#3#4
{\biggl\{ \hs{#1}\setbar
{\matrix{\hbox{#2} \hfill\cr\hbox{#3} \hfill\cr\hbox{#4} \hfill}}
\biggr\} }

\def\locus#1{\{ #1 \}}

\def\varvarmatrix#1#2#3#4
{\def\normalbaselines{\baselineskip#1\lineskip3pt\lineskiplimit3pt}
\vbox{\vskip #2 \hbox{\hfill$\matrix{#3}$\hfill}\vskip #3}}

\def\varmatrix#1#2
{\def\normalbaselines{\baselineskip#1\lineskip3pt\lineskiplimit3pt}
\matrix{#2}}

\def\bigcases#1
{\biggl\{\, \vcenter{\normalbaselines{\mathsurround=0pt}
\ialign{$##\hfil$&\quad##\hfil\crcr#1\crcr}}\biggr.}

\def\st{\subset}
\def\sm{\setminus}

\def\pione{\pi\sb1}

\def\vol#1{{\bf #1 }}

\def\and{\hbox{and}}
\def\where{\hbox{where}}

\def\afcite{\ }

\def\Pt{\P\sp 2}
\def\hc{X, Y, Z}
\def\F{\Scr F}
\def\L{\Scr{L}}

\def\W{\Scr{W}}

\def\pqi{\phi\sb q \inv}

\def\ArtalBartolo{[1]}
\def\Libgober{[2]}
\def\Nemethi{[3]}
\def\OkaS{[4]}
\def\OkaPre{[5]}
\def\ShimadaF{[6]}
\def\ShimadaW{[7]}
\def\Tokunaga{[8]}
\def\ZariskiE{[9]}
\def\ZariskiP{[10]}
\def\ZariskiD{[11]}

\def\ShimadaFPart#1{[6; #1]}
\def\ShimadaWPart#1{[7; #1]}

\def\sbinfty{\sb{\infty}}

\def\spno{\sp{n-1}}

%
%
%

\centerline{\bf A note on Zariski pairs}
\bigskip
\centerline{Ichiro Shimada}
\smallskip
\centerline{Max-Planck-Institut f\"ur Mathematik, Bonn}

\bign
{\bf \S 0. Introduction}
\medskip
In \ArtalBartolo,
Artal Bartolo defined the notion of {\sl Zariski pairs} as follows:
\medskip
{\bf Definition.}
\hs
A couple of complex reduced projective
plane curves $C\sb 1$ and $C\sb 2$ of a same degree
is said to make a  Zariski pair,
if there exist tubular neighborhoods
$T(C\sb i) \st \Pt $ of $C\sb i$ for $i=1,2$
such that $(T(C\sb 1), C\sb 1)$ and $(T(C\sb 2 ), C\sb 2)$
are diffeomorphic,
while the pairs $(\Pt, C\sb 1)$ and $(\Pt, C\sb 2)$ are  not homeomorphic;
that is, the singularities of $C\sb 1 $ and $C\sb 2$
are topologically equivalent,
but the embeddings of $C\sb 1$ and $C\sb 2$
into $\Pt$ are not topologically equivalent.
\medskip
The first example of Zariski pair was discovered and studied
by Zariski in \ZariskiE\afcite and \ZariskiD.
He showed that
there exist projective plane curves $C\sb 1 $ and $C\sb 2$ of degree $6$
with $6$ cusps and no other singularities
such that $\pione (\Pt \sm C\sb 1)$ and $\pione (\Pt \sm C\sb 2)$
are not isomorphic.
Indeed, the placement of the $6$ cusps  on the sextic curve  has a
crucial effect on
the fundamental group of the complement.
Let $C\sb 1$ be a sextic curve defined by an equation $f\sp 2 + g \sp 3= 0$,
where $f$ and $g$ are general homogeneous polynomials of degree $3$ and $2$,
respectively.
Then $C\sb 1$ has $6$ cusps lying on a conic defined by $g=0$.
In \ZariskiE, it was shown
that  $\pione (\Pt\sm C\sb 1)$ is isomorphic to the free product
$\Z/(2) * \Z/(3)$
of cyclic groups of order $2$ and $3$.
On the other hand,
in \ZariskiD, it was proved that there exists  a sextic curve $C\sb 2$
with $6$ cusps which are not  lying on any conic, and
that the fundamental group $\pione (\Pt \sm C\sb 2)$ is cyclic of order $6$.
In \OkaS, Oka gave an explicit defining equation of $C\sb 2$.
In \ArtalBartolo,
Artal Bartolo presented  a simple way to construct
$(C\sb 1, C\sb 2)$
from a cubic curve $C$
by means of a Kummer covering of $\P\sp 2$
of exponent $2$
branched along three lines tangent to $C$ at its points of inflection.
\medskip
Except for this example,
very few Zariski pairs are known (\ArtalBartolo, \Tokunaga).
In \OkaPre, and independently in \ShimadaW,
infinite series of Zariski pairs
have been constructed from the above example of Zariski by means of
covering tricks of the plane.
\medskip
In this paper,
we present a method to construct Zariski pairs,
which yields
two infinite series of new examples of Zariski pairs as special cases.
\medskip
A germ of curve singularity is called of type $(p, q)$
if it is locally defined by $x\sp p + y\sp q=0$.
\medskip
{\it Series I.}
\hs
This series consists of pairs $(C\sb 1 (q), C\sb 2 (q))$
of curves of degree $3q$, where
$q$ runs through the set of integers $\ge 2$
prime to $3$.
Each of $C\sb 1 (q)$ and $C\sb 2 (q)$
has $3q$ singular points of type $(3,q)$
and no other singularities.
The fundamental group $\pione (\P\sp 2 \sm C\sb 1 (q))$
is non-abelian,
while
$\pione (\P\sp 2 \sm C\sb 2 (q))$
is abelian.
When $q=2$,
this example is nothing but  the classical one
of the sextic curves due to Zariski.
\medskip
{\it Series II.} \hs
This series consists of pairs $(D\sb 1 (q), D\sb 2 (q))$
of curves of degree $4q$, where
$q$ runs through the set of  odd integers $> 2$.
Each of $D\sb 1 (q)$ and $D\sb 2 (q)$
has $8q$ singular points of type $(2,q)$ -- that is,
rational double  points of type $A\sb{q-1}$ --
and no other singularities.
The fundamental group $\pione (\P\sp 2 \sm D\sb 1 (q))$
is non-abelian,
while
$\pione (\P\sp 2 \sm D\sb 2 (q))$
is abelian.
\medskip
Our method is a generalization of Artal Bartolo's method
for re-constructiong the classical example of Zariski
to higher dimensions and arbitrary exponents of the Kummer covering.
Indeed, when $q=2$ in Series I,
our construction coincides with his.
\medskip
Instead of  the computation
of the first Betti number
of the cyclic branched covering
of $\P\sp 2$,
which was employed in \ArtalBartolo,
we use the fundamental groups of the complements
in order to distinguish two embeddings of curves in $\Pt$.
For the calculation of the fundamental groups,
we use \ShimadaFPart{Theorem 1}
and a result of \Nemethi\afcite and \ShimadaW.
\medskip
{\sl Acknowledgment.}
\hs
Part of this work was done during the author's stay
at Institute of Mathematics in Hanoi
and Max-Planck-Institut f\"ur Mathematik  in Bonn.
The author thanks to people at these institutes
for their warm hospitality.
He also thanks
to Professors M.\ Oka  and H.\ Tokunaga
for stimulating discussions.
\bign
{\bf \S 1. A method of constructing Zariski pairs}
\medskip
{\bf  1.1. Non-abelian members.}
\hs
Let $p$ and $q$ be integers $\ge 2$
prime to each other.
We choose homogeneous polynomials
$f \in H\sp 0 (\Pt, \O(pk))$ and $g \in H\sp 0 (\Pt, \O (qk))$,
where $k$ is an integer $\ge 1$.
Suppose that $f$ and $g$ are generally chosen.
Consider the projective plane curve
$$
C\sb{p,q, k} \quad : \quad f\sp q + g\sp p = 0
$$
of degree $pqk$ (cf.\hs  \Libgober).
It is easy to see that the singular locus of this curve
consists of $pqk\sp 2 $ points of type $(p, q)$.
In \ShimadaWPart{Example (3) in \S 0},
the following is shown.
\medn
{\bf Proposition 1.}
\hs
{\sl
The fundamental group $\pione (\Pt\sm C\sb{p, q, k})$
is isomorphic to the group
$\langle \hs a, b, c \hs | \hs a\sp p = b\sp q = c,
\hs c\sp k = 1 \hs \rangle$.
In particular, it is non-abelian.
}
\medn
See also \Nemethi, in which the fundamental groups
of the complements of curves of this type are  calculated.
There the groups are presented in a different way.
\medskip
This curve $C\sb{p, q, k}$ will be a member $C\sb 1$ of a Zariski pair.
\medskip
{\bf 1.2. Abelian partners.}
\hs
We shall construct the other  member $C\sb 2 $
of the Zariski pair such that
$\pione (\Pt \sm C\sb 2)$ is abelian.
\medskip
Let $p$, $q$ and $k$ be integers as above.
We put $n = pk$.
Interchanging $p$ and $q$ if necessary,
we may assume that $n \ge 3$.
Let $S\sb 0 \st \P\sp{n-1}$ be a hypersurface of degree $n$
defined by
$F\sb 0 (X\sb 1, \dots, X\sb n) = 0$.
We consider a linear pencil of hypersurfaces
$$
S\sb t \quad : \quad F\sb 0 (X\sb 1, \dots, X\sb n) +
t \cdot X\sb 1 \cdots X\sb n =0,
$$
which is spanned by $S\sb 0$ and $S\sb{\infty} :=
\locus {X\sb 1 \cdots X\sb n =0}$.
We put $H\sb i = \locus{X\sb i = 0}$ $(i= 1, \dots, n)$.
We consider the morphism
$\phi\sb q : \P\sp {n-1} \to \P\sp{n-1}$
given by
$$
(Y\sb 1 : \dots : Y\sb n) \quad \longmapsto
\quad ( X\sb 1 : \dots : X\sb n ) \hs =\hs
( Y\sb 1 \sp q : \dots : Y\sb n \sp q),
$$
which is a covering of degree $q\spno$ branched along $S\sbinfty$.
\medn
{\bf Proposition 2.}
\hs
{\sl
Suppose that (1) every member $S\sb t$
is reduced,
and that (2) $S\sb 0$ contains none of the hyperplanes $H\sb i$.
Then
$\pione (\P\spno\sm \phi\sb q \inv (S\sb t))$ is abelian
for a general member $S\sb t$.
}
\medn
{\it Proof.}
\hs
Let $\P\sp 1$ be the $t$-line,
and we put $\A\sp 1 := \P\sp 1 \sm \locus {\infty}$.
Let $\W \st \P\spno \times \A\sp 1$
be the divisor defined by
$$
X\sb 1 \cdots X\sb n \cdot
(\hs  F\sb 0 (X\sb 1 , \dots , X\sb n) + t\cdot X\sb 1 \cdots X\sb n \hs ) = 0,
$$
which is the union of $S\sbinfty \times \A\sp 1 $
and the universal family of the
affine part $\set{ S\sb t }{t \in \A\sp 1}$
of the pencil.
For $t \in \A\sp 1$,
we denote by $W\sb t \st \P\spno$
the divisor obtained from the scheme theoretic intersection
$ (\locus{t} \times \P\spno ) \cap \W$,
which is equal with the divisor $S\sb t + S\sbinfty$.
\smallskip
First, we shall show that $\pione (\P\spno \sm W\sb t)$
is abelian for a general $t$.
Remark  that the assumption (2) implies that $S\sb t$ contains
none of $H\sb i$ unless $t=\infty$.
Combining this with the assumption (1),
we see that $W\sb t$ is reduced for all $t\in \A\sp 1$.
Hence, by \ShimadaFPart{Theorem 1},
the inclusion
$\P\spno \sm W\sb t  \hookrightarrow (\P\spno \times \A\sp 1 ) \sm \W$
induces an isomorphism on the fundamental groups
for a general $t$.
Therefore,
it is enough to show that
$\pione ((\P\spno \times \A\sp 1 )\sm \W)$
is abelian.
In order to prove this,
we consider the first projection
$$
p \quad :\quad
(\P\spno \times \A\sp 1 )\sm \W \quad \lra \quad \P\spno \sm S\sbinfty.
$$
Since $\set{ S\sb t }{t \in \P\sp 1}$
is a pencil whose base locus is contained in $S\sbinfty$,
there is a unique member $S\sb{t(P)}$ $(t(P)\not = \infty)$
containing $P$
for each point $P\in \P\spno \sm S\sbinfty$.
Therefore
$p\inv (P)$ is a punctured affine line $\A\sp 1 \sm\locus {t(P)}$
for every $P\in \P\spno \sm S\sbinfty$.
Consequently,  $p$ is a locally trivial fiber space.
Moreover,
$p$ has a section
$$
s \quad : \quad
\P\spno \sm S\sbinfty \quad\lra\quad (\P\spno \times \A\sp 1 )\sm \W,
$$
which is given by, for example,
$s(P) = (P, \hs t(P) +1)$.
Hence the homotopy exact sequence of $p$ splits.
Combining this with the fact that the image of the injection
$\pione (\A\sp 1 \sm \{ t(P) \}) \to \pione ((\P\spno \times \A\sp 1 )\sm\W)$
is contained in the center, we see that
$$
\pione ((\P\spno \times \A\sp 1 )\sm \W)
\quad \cong\quad
\pione (\P\spno \sm S\sbinfty )\hs
\times \hs \pione (\A\sp 1 \sm \locus {\hbox { a point }}).
$$
This shows that
$\pione ((\P\spno \times \A\sp 1 )\sm \W)$
is abelian.
\smallskip
Note that
$\phi\sb q : \P\spno \to \P\spno$
is \'etale over $\P\spno \sm W\sb t$ for every $t$.
Hence the natural homomorphism
$$
\phi\sb{q*} \quad :
\quad \pione (\P\spno \sm \phi\sb q \inv (W\sb t)) \quad \lra \quad
\pione (\P\spno \sm W\sb t)
$$
is injective.
This implies that
$\pione (\P\spno\sm \phi\sb q \inv (W\sb t))$
is abelian for a general $t$.
On the other hand,
since
$\P\spno \sm \phi\sb q \inv (W\sb t)$
is a Zariski open dense subset of $\P\spno \sm \phi\sb q \inv (S\sb t)$,
the inclusion
induces a surjective homomorphism
$$
\pione (\P\spno \sm \phi\sb q \inv (W\sb t)) \quad\surj\quad
\pione (\P\spno \sm \phi\sb q \inv (S\sb t)).
$$
Thus  $\pione (\P\spno \sm \phi\sb q \inv (S\sb t))$ is also
abelian for a general $t$.
\qed
\medn
{\bf Proposition 3.}
\hs
{\sl
Suppose the following;
(3) $S\sb 0 \cap H\sb i$
is a non-reduced divisor $p D\sb i$ of $H\sb i$ of multiplicity $p$,
where $D\sb i$ is a reduced divisor of $H\sb i$,
none of  whose irreducible components
is contained in $H\sb i \cap (\cup\sb{j\not = i} H\sb j)$,
and
(4) the singular locus of $S\sb t$ is of codimension $\ge 2$ in $S\sb t$
for a general $t$.
Then the general plane section
$\P\sp 2 \cap \phi \sb q \inv (S\sb t)$ of $\pqi (S\sb t)$
is a curve of degree $pqk$,
and its singular locus consists of $pqk\sp 2 $ points of type $(p, q)$.
}
\medn
{\it Proof.}
\hs
Note that the assumption (3)
implies that $S\sb t \cap H\sb i$
is also equal with $p D\sb i$ for $t\not = \infty$.
Let $P$ be a general point of any irreducible component of $D\sb i$,
and let $Q$ be a point such that $\phi\sb q (Q) = P$,
which lies on the hyperplane defined by $Y\sb i = 0$.
By the assumption (3),
$Q$ is not contained in any of the other hyperplanes
defined by $Y\sb j = 0$ $(j\not = i)$.
Hence there exist analytic local coordinate systems
$(w\sb 1 , \dots, w\sb{n-1})$
and $(z\sb 1 , \dots, z\sb{n-1})$
of $\P\spno$ with the origins $P$ and $Q$,
respectively,
such that
$H\sb i$ is given by $w\sb 1 =0$,
$\phi\sb q \inv (H\sb i)$
is given by $z\sb 1 =0$,
and $\phi\sb q$
is given by
$$
(z\sb 1 , \dots, z\sb{n-1} ) \quad \longmapsto
\quad (w\sb 1, \dots, w\sb{n-1}) \hs =
\hs  (z\sb 1 \sp q , z\sb 2 , \dots, z\sb{n-1}).
$$
Let $t \in \A\sp 1$ be general.
By the assumption (3),
the defining equation of $S\sb t$ at $P$ is of the form
$$
u(w)\cdot w\sb 1 + v(w\sb 2, \dots, w\sb{n-1}) \sp p = 0.
$$
By the assumption (4),
$S\sb t$ is non-singular at $P$,
because $P$ is a general point
of an irreducible component of $D\sb i$.
This implies that
$u(P) \not= 0$.
On the other hand,
the divisor $D\sb i$,
which is defined by $v(w\sb 2 , \dots, w\sb{n-1})=0$
on the hyperplane $H\sb i = \locus{w\sb 1 =0}$,
is non-singular at $P$,
because $D\sb i$ is reduced
by the assumption (3) and $P$ is general.
Hence we have
$$
{{\partial v}\over{\partial w\sb j}} (P) \not= 0
\quad\hbox{ at least for one $j \ge 2$.}
$$
The defining equation of $\phi\sb q\inv (S\sb t)$
is then of the form
$$
\wt u (z) \cdot z\sb 1 \sp q + v (z\sb 2 , \dots, z\sb{n-1} )\sp p = 0,
\quad\where\quad \wt u (Q) \not = 0.
$$
Then, it is easy to see that,
in terms of suitable analytic coordinates $(\wt z\sb 1 , \dots,
\wt z\sb{n-1})$
with the origin $Q$,
this equation can be written
as follows;
$$
\wt z\sb 1 \sp{\hskip 2pt q} \hs + \hs  \wt z \sb 2 \sp{\hskip 2pt p}
\hs = \hs 0.
$$
Thus, when we cut $\phi\sb q \inv (S\sb t)$
by a general $2$-dimensional plane passing through $Q$,
a germ of curve singularity of type $(p, q)$ appears at $Q$.
\smallskip
Since the degree of $D\sb i$ is $k = n/p$,
the inverse image $\phi\sb q \inv (D\sb i)$
is a reduced hypersurface of degree $qk$ in
the hyperplane defined by $Y\sb i=0$.
Moreover $\pqi (D\sb i)$ and $\pqi (D\sb j)$
have no common irreducible components
when $i\not = j$
because of the assumption (3).
Hence  the intersection  points of
$\phi\sb q \inv (\sum \sb{i=1} \sp n D\sb i)$
with a general plane $\Pt \st \P\sp n$
is $pqk\sp 2$ in number.
Moreover, $\Pt \cap \phi\sb q \inv (S\sb t)$
is non-singular outside of these intersection points,
because of the assumption (4).
\qed
\medskip
{\bf  1.3. Summary.}
\hs
Suppose that we have constructed a hypersurface $S\sb 0 \st \P\spno$
of degree $n\ge 3$ which satisfies the assumptions (1)-(4)
in Propositions 2  and 3.
Let $C\sb 2$ be a general plane section of $\pqi (S\sb t)$,
where $t$ is general.
Because of Zariski's hyperplane section theorem
\ZariskiP\afcite  and Propositions 1, 2 and 3,
we see that
the curve $C\sb 2$ has the same type of singularities as that of
$C\sb{ p, q, k}$,
but the fundamental group $\pione (\Pt \sm C\sb 2)$ is abelian.
Hence $(C\sb 1, C\sb 2)$ is a Zariski pair, with $C\sb 1 = C\sb{p, q, k}$.
\bign
{\bf \S 2.  Construction of Series I}
\medskip
We carry out the construction of the previous section with
$p=3$, $k=1$, $n=3$  and
$q$ an arbitrary integer $\ge 2$ prime to $3$.
\medskip
We fix a homogeneous coordinate system $(X:Y:Z)$
of $\Pt$, and put
$$
\varmatrix{20pt}{
&L\sb1  = \locus{X=0}, \quad L\sb2  = \locus{Y=0},
\quad L\sb 3  = \locus {Z=0}, \quad \and \cr
&R\sb 1 = (0:1:-1) \hs \in \hs L\sb 1,
\quad
R\sb 2 = (1:0:-1) \hs \in \hs L\sb 2.
}
$$
Let $\P\sb * (\GG (\Pt, \O(3)))$ be the space of all cubic curves on $\Pt$,
which is isomorphic to the projective space of dimension $9$,
and  let $\F \st \P\sb * (\GG (\Pt, \O(3)))$ be the family of
cubic curves $C$
which satisfy the following conditions;
\smalln
(a) $C$ intersects $L\sb 1$ at $R\sb 1$ with multiplicity $ \ge 3$,
\parn
(b) $C$ intersects $L\sb 2$ at $R\sb 2$ with multiplicity $\ge 3$, and
\parn
(c) $C$ intersects $L\sb 3$ at a point  with multiplicity $\ge 3$.
\smalln
(We consider that $C$ intersects a line $L\sb i$
with multiplicity $\infty$, if $L\sb i $ is contained in $C$.)
\medn
{\bf Proposition 4.}
\hs
{\sl
The family $\F$ consists of $3$ projective lines.
They meet at one point corresponding to
$C\sbinfty : = \locus {XYZ=0}$.
}
\medn
{\it Proof.}
\hs
Let $F(\hc)=0$ be the defining equation of a member $C$ of this family $\F$.
By the condition (a),
$F$ is of the form
$$
F(\hc) = A (Y+Z) \sp 3 + X \cdot G(\hc),
$$
where $A$ is a constant,
and $G(\hc)$ is a homogeneous polynomial of degree $2$.
By the condition (b), we have
$F(X, 0, Z) = A (Z+X)\sp 3$,
and hence we get
$$
G(\hc) = A (3Z\sp 2 + 3ZX + X\sp 2 ) + Y  \cdot H(\hc),
$$
where $H (\hc) $ is a homogeneous polynomial of degree $1$.
By the condition (c), we have
$F(X, Y, 0) = A (Y+ \aa X)\sp 3$
for some $\aa$. Then $\aa$ must be a cubic root of unity, and we get
$$
H(\hc) = 3A\aa \sp 2 X + 3A \aa Y + BZ,
$$
where $B$ is a constant.
Combining all of these, we get
$$
\eqalign{
&F(\hc) \cr
= &A(X\sp 3 + Y\sp 3 + Z\sp 3) + 3 A ( \aa\sp 2 X\sp 2 Y
+ \aa X Y \sp 2 +  Y\sp 2 Z +  Y Z \sp 2 +  Z\sp 2 X +  Z X \sp 2)
+ BXYZ\cr
=& A(X+Y+Z) \sp 3 +
3 A (\aa\sp 2 -1 ) X\sp 2 Y + 3 A (\aa -1) X Y\sp 2 + (B-6A) XYZ.
}
$$
This  curve  $C = \locus{F=0}$ intersects $L\sb 3$ at
$$
R\sb 3 \hs = \hs  R\sb 3 (\aa) \hs := \hs (1: -\aa: 0) \hs \in \hs L\sb 3
$$
with multiplicity $\ge 3$.
This means that the family $\F$
consists of three lines $\L (1)$, $\L(\ww)$ and $\L(\ww\sp 2)$
in the projective space  $\P\sb * (\GG (\Pt, \O(3)))$,
where $\ww= \exp(2\pi i / 3)$,
such that a general cubic
$C$ in $\L (\aa)$ intersects $L\sb 3$ at $R\sb 3 (\aa)$
with multiplicity $3$.
The ratio of the coefficients $t:= B/ A$
gives an affine coordinate on each line $\L (\aa)$.
The three lines $\L (1)$, $\L(\ww)$, $\L(\ww\sp 2)$
intersect at one point $t = \infty$
corresponding to the cubic $C\sb{\infty}=L\sb 1 + L\sb2 + L\sb 3$.
\qed
\medskip
Hence we get three pencils of cubic curves
$\set{ C(1)\sb t }{t \in \L (1)}$,
$\set{ C(\ww)\sb t }{t \in \L (\ww)}$, and
$\set{ C(\ww\sp 2 )\sb t }{t \in \L (\ww\sp 2)}$.
It is easy to check that these pencils
satisfy
the assumptions (2), (3) and (4)
in the previous section.
Note that  the pencil $\L (1)$
does not satisfy the assumption (1)
because $C(1)\sb 6 $ is a triple line.
However, the other two satisfy (1).
Indeed, if a cubic curve
$C$ in the family $\F$ is non-reduced,
then the conditions (a)-(c) imply
that it must be a triple line. Therefore the three points
$R\sb 1$, $ R\sb 2$ and $ R\sb 3 (\aa) $ are co-linear,
which is equivalent to $\aa=1$.
Consequently, $C$ must be a member of $\L (1)$.
\medskip
Now,
by using the pencil $\L (\ww)$
or $\L(\ww\sp 2)$, we complete  the construction of Series I.
\medskip
Note that,
if $C(1)\sb a$ is a non-singular member of $\L (1)$,
then
$\pione (\Pt \sm \phi\sb q \inv (C(1)\sb a))$
is isomorphic to the free product $\Z / (3) * \Z / (q)$.
Indeed, since
$C(1)\sb a$ is defined by
$$
(X+Y+Z)\sp 3 + (a-6)XYZ = 0,
$$
the pull-back
$\phi\sb q\inv (C(1)\sb a)$
is defined by
$$
(U\sp q + V\sp q + W\sp q)\sp 3 + (a-6) (UVW)\sp q =0,
$$
which is of the form $\wt f\sp{\hskip 2pt 3} + \wt g\sp{\hskip 2pt q} =0$.
The polynomials $\wt f$ and $\wt g$ are
not general by any means.
However,
since the type of singularities of $\phi\sb q \inv (C(1)\sb a )$
is the same as that of $C\sb{3,q,1}$,
we have an isomorphism
$\pione (\Pt\sm \phi\sb q \inv (C(1)\sb a)) \cong \pione
(\Pt \sm C\sb{3,q,1})$.
\bign
{\bf \S 3. Construction of Series II}
\medskip
It is enough to show the following:
\medn
{\bf Proposition 5.}
\hs
{\sl
The quartic surface
$$
S\sb 0 \quad:\quad
F\sb 0 (x\sb 1, x\sb 2, x\sb 3, x\sb 4) \hs :
= \hs (x\sb 1 \sp 2 + x\sb 2 \sp 2 )\sp 2 +
2 x\sb 3 x\sb 4 (x\sb 1 \sp 2 - x\sb 2 \sp 2 ) + x\sb 3 \sp 2 x\sb 4 \sp 2
\hs = \hs 0
$$
in $\P\sp 3 $ satisfies the assumptions (1)-(4) with $p=2$ and $k=2$.
}
\medn
{\it Proof.}
\hs
The assumptions (2) and (3) can be trivially checked.
To check the assumptions (1) and (4),
we put
$$
F\sb t  \hs  := \hs F\sb 0 + t\cdot x\sb 1 x\sb 2 x\sb 3 x\sb 4,
$$
and calculate the partial derivatives $\partial F\sb t / \partial x\sb i$
for $i=1, \dots, 4$.
Let $Q\sb t \st \P\sp 3$
be the quadric surface defined by
$$
 2 x\sb 1 \sp 2 - 2 x\sb 2 \sp 2 + 2 x\sb 3 x\sb 4 + t x\sb 1 x\sb 2 \hs
= \hs 0.
$$
It is easy to see that $Q\sb t$ is irreducible for all $t\not = \infty$.
It is also easy to see that $Q\sb t$ is the unique common
irreducible component
of the two cubic surfaces
$$
{{\partial F\sb t }\over{\partial x\sb 3 }} = 0,
 \qquad\and\qquad
{{\partial F\sb t }\over{\partial x\sb 4 }} = 0.
$$
Suppose that a surface $S\sb a = \locus{F\sb a =0}$
in this pencil contains a non-reduced
irreducible component $m T$ $(m \ge 2)$.
Then, both of $\partial F\sb a / \partial x\sb 3$ and
$\partial F\sb a / \partial x\sb 4$
must vanish on $T$.
Hence $T$ must coincide with $Q\sb a $, and we get $S\sb a = 2 Q\sb a$.
Comparing the defining equations of $S\sb a $ and $2 Q\sb a$,
we see that there are no such $a$.
Thus the assumption (1) is satisfied.
To check the assumption (4),
we remark that the condition $\dim\Sing S\sb t \le 0$
is an open condition for $t$.
Hence it is enough to prove, for example,
$\dim \Sing S\sb 2 =0$.
It is easy to show that $\Sing S\sb 2$
consists of four points
$(1:\pm\sqrt{-1}:0:0)$, $(0:0:0:1)$ and $(0:0:1:0)$.
\qed
\bigskip
\centerline{\bf References}
\medskip
\item{\ArtalBartolo}
 Artal Bartolo, E.:
{\sl Sur les couples de Zariski},
J.\ Alg.\ Geom. {\bf 3} (1994), 223 - 247
\item{\Libgober}
Libgober, A.:
{\sl Fundamental groups of the complements to plane singular curves},
Proc.\ Symp.\ in Pure Math.
{\bf 46} (1987), 29 - 45
\item{\Nemethi}
N\'emethi, A.:
{\sl On the fundamental group of the complement
of certain singular plane curves},
Math. Proc. Cambridge Philos. Soc.
{\bf 102} (1987), 453 - 457
\item{\OkaS}
Oka, M.:
{\sl Symmetric plane curves with nodes and cusps},
J.\ Math.\ Soc.\ Japan {\bf 44} (1992), 375 - 414
\item{\OkaPre}
Oka, M.: {\sl Two transformations of plane
curves and their fundamental groups},
preprint
\item{\ShimadaF}
Shimada, I.:
{\sl Fundamental groups of open algebraic varieties},
to appear in Topology
\item{\ShimadaW}
Shimada, I.:
{\sl A weighted version of Zariski's hyperplane section theorem
and fundamental groups of complements of plane curves},
preprint
\item{\Tokunaga}
Tokunaga, H.:
{\sl A remark on Bartolo's paper},
preprint
\item{\ZariskiE}
Zariski, O.:
{\sl On the problem of existence of
algebraic functions of two variables
possessing a given branch curve},
Amer. J. Math. {\bf 51} (1929), 305 - 328
\item{\ZariskiP}
Zariski, O.:
{\sl A theorem on the Poincar\'e group of an algebraic hypersurface},
Ann.\ Math. {\bf 38} (1937), 131 - 141
\item{\ZariskiD}
Zariski, O.:
{\sl The topological discriminant group of a Riemann surface
of genus $p$},
Amer.\ J.\ Math.\ \vol{59} (1937), 335 - 358
\bign
Max-Planck-Institut f\"ur Mathematik
\parn
Gottfried-Claren-Strasse 26
\parn
53225 Bonn, Germany
\parn
shimada@mpim-bonn.mpg.de

\end